\documentclass{article}
\usepackage{spconf,amsmath,amssymb,graphicx,xcolor}
\usepackage[linesnumbered,ruled,vlined]{algorithm2e}
\usepackage{algpseudocode}
\usepackage{enumerate}
\usepackage{xfrac}

\def\L{{\boldsymbol{\Lambda}}}

\def\cE{\mathcal{E}}

\def\cG{\mathcal{G}}

\def\cL{\mathcal{L}}

\def\cN{\mathcal{N}}
\def\cO{\mathcal{O}}

\def\smskip{\smallskip}

\def\texitem#1{\par\smskip\noindent\hangindent 25pt
               \hbox to 25pt {\hss #1 ~}\ignorespaces}

\def\norm#1{\left\|#1\right\|}

\newcommand{\BEAS}{\begin{eqnarray*}}
\newcommand{\EEAS}{\end{eqnarray*}}
\newcommand{\BEA}{\begin{eqnarray}}
\newcommand{\EEA}{\end{eqnarray}}
\newcommand{\BEQ}{\begin{eqnarray}}
\newcommand{\EEQ}{\end{eqnarray}}
\newcommand{\BIT}{\begin{itemize}}
\newcommand{\EIT}{\end{itemize}}
\newcommand{\BNUM}{\begin{enumerate}}
\newcommand{\ENUM}{\end{enumerate}}

\newcommand{\BA}{\begin{array}}
\newcommand{\EA}{\end{array}}

\newcommand{\ones}{\mathbf 1}

\newcommand{\reals}{\mathbb{R}}
\newcommand{\integers}{\mathbb{Z}}

\newcommand{\Rank}{\mathop{\bf rank}}

\newcommand{\diag}{\mathop{\bf diag}}

\newif\ifpagenumbering
\pagenumberingtrue

\pagenumberingfalse

\newsavebox{\theorembox}
\newsavebox{\lemmabox}
\newsavebox{\defnbox}
\newsavebox{\corollarybox}
\newsavebox{\propositionbox}
\newsavebox{\remarkbox}
\newsavebox{\assbox}
\savebox{\theorembox}{\noindent\bf Theorem}
\savebox{\lemmabox}{\noindent\bf Lemma}
\savebox{\defnbox}{\noindent\bf Definition}
\savebox{\corollarybox}{\noindent\bf Corollary}
\savebox{\propositionbox}{\noindent\bf Proposition}
\savebox{\remarkbox}{\noindent\bf Remark}
\savebox{\assbox}{\noindent\bf Assumption}

\DeclareMathOperator*{\argmin}{\arg\!\min}

\def\id{\mathbf{I}}
\def\zero{\mathbf{0}}
\def\one{\mathbf{1}}

\def\sym{\mathbb{S}}

\def\L{{\cal L}}
\def\N{\mathbf{Null}}

\title{\vspace*{-8mm}Decentralized Computation of Effective Resistances\\ and Acceleration of Consensus Algorithms}
\name{N.S. Aybat   M. G\"urb\"uzbalaban*}
\address{Author Affiliation(s)}
\twoauthors
  {Necdet Serhat Aybat\sthanks{Research of N. S. Aybat was partially supported by NSF grants CMMI-1400217 and CMMI-1635106, and ARO grant W911NF-17-1-0298.}}
	{Department of IME\\
	Pennsylvania State University\\
	University Park, PA, USA\\
	nsa10@psu.edu}
  {Mert G\"urb\"uzbalaban\sthanks{Research of M. G\"urb\"uzbalaban was partially supported by the NSF grant DMS-1723085.}}
	{Department of MSIS\\
	Rutgers Business School\\
	Piscataway, NJ, USA\\
	mg1366@rutgers.edu}
\begin{document}
%
\maketitle
\begin{abstract}
The effective resistance between a pair of nodes {in} a weighted undirected graph is defined as the 
potential difference induced between them when a unit current is injected at the first node and extracted at the second node, {treating edge weights as the conductance values of edges}. The effective resistance is a key quantity of interest in many applications and fields including solving linear systems, Markov Chains and continuous-time averaging networks. We develop an efficient linearly convergent distributed algorithm for computing effective resistances and demonstrate its performance through numerical studies. We also apply our algorithm to the consensus problem where the aim is to compute the average of node values in a distributed manner. We show that the distributed algorithm we developed for effective resistances can be used to 
accelerate the convergence of the classical consensus iterations considerably by a factor depending on the network structure.
\end{abstract}
\begin{keywords}
Effective resistance, graph, distributed optimization, consensus, Laplacian matrix, Kaczmarz method
\vspace*{-3mm}
\end{keywords}
\vspace*{-3mm}
\section{Introduction}\label{sec:intro}
Let $\cG=(\mathcal{N},\mathcal{E},w)$ be an undirected, weighted and connected graph defined by the set of nodes (agents) $\cN=\{1,\ldots,n\}$, the set of edges $\cE\subset\cN\times\cN$, and 
the edge weights $w_{ij} > 0$ for $(i,j) \in \cE $. Since $\cG$ is undirected, we assume that both $(i, j)$ and $(j, i)$ refer to the same edge when it exists, and for all $(i,j)\in\mathcal{E}$, we set $w_{ji} = w_{ij}$. Identifying the weighted graph $\cG$ as an electrical network 
in which each edge $(i,j)$ corresponds to a branch of conductance $w_{ij}$, the effective resistance $R_{ij}$ between a pair of nodes $i$ and $j$ is defined as the voltage potential difference induced between them when a unit current is injected at $i$ and extracted at $j$.

The effective resistance, also known as the resistance distance, is a key quantity of interest to compute in many applications and algoritmic questions over graphs. It defines a metric on graphs providing bounds on its conductance \cite{klein2002resistance,Klein1993}. Furthermore, it is closely associated with the hitting time and commute time for a random walk\footnote{The hitting time $H_{ij}$ is the
expected number of steps of a random walk starting from 
$i$ until it first
visit 
$j$. The commute time $C_{ij}$ is the expected
number of steps required to go from $i$ to $j$ and back again.}
on the graph $G$ such that the probability of a transition from $i$ to $j^*\in\cN_i$ is $w_{ij^*}/\sum_{j\in\cN_i}w_{ij}$ where {$\cN_i\triangleq\{j\in\cN:~(i,j)\in\cE\}$} 
denotes the set of neighboring nodes of $i\in\cN$; therefore, it arises naturally for  studying random walks over graphs and their mixing time properties \cite{boydSIAMreview,aldous-fill-2014,doyle1984random}, continuous-time averaging networks including consensus problems in distributed optimization \cite{boydSIAMreview}. Other prominent applications include distributed control and estimation \cite{BarooahCDC}, solving symmetric diagonally dominant (SDD) linear systems \cite{SpielmanSrivastava}, deriving complexity bounds in the Asymmetric Traveling Salesman Problem (ATSP) \cite{AnariTSP}, design and control of communication networks \cite{Tizghadam10,jadbabaie2004geographic} and spectral sparsification of graphs \cite{kapralov2012spectral}. 

There exist \emph{centralized} algorithms for computing or approximating 
{$\{R_{ij}\}_{i\neq j}$} accurately which require global communication beyond local communication among the neighboring agents \cite{SpielmanSrivastava,bapat2003simple}. They are based on computing or approximating the entries of the pseudoinverse $\cL^{+}$ of the Laplacian matrix, based on the identity $ R_{ij} = \cL^{+}_{ii} + \cL^{+}_{jj} - 2 \cL^{+}_{ij} $
\cite{SpielmanSrivastava}. However, such centralized algorithms are \emph{impractical} or \emph{infeasible} for several key applications in multi-agent systems where only local communications between the neighboring agents are allowed; this motivates the development of distributed algorithms for computing effective resistances, {which are used in solving many optimization and estimation problems over graphs.}
Prominent examples include, least square regression and more general estimation problems over graphs, formation control of moving agents with noisy measurements and stability of multi-vehicle swarms \cite{BarooahCDC}.

To our knowledge, there has been no systematic study of distributed algorithms for computing effective resistances. In this work, we discuss how existing algorithms in the distributed optimization literature for solving linear systems can be adapted to solve this problem. First, we show that a naive implementation of {consensus optimization methods, e.g.,} the EXTRA algorithm~\cite{shi2015extra} 
is inefficient in terms of the convergence and communication requirements. Second, we propose a variant of the Kaczmarz method and show that it is linearly convergent while being efficient in terms of total number of local communications carried out. Third, we demonstrate the performance of our algorithms on numerical examples. In particular, numerical experiments suggest finite convergence of our algorithms which is of independent interest. {Finally, we apply our results to the consensus problem \cite{boyd2006randomized} where the aim is to compute the average of values assigned to each node in a distributed manner. 
Specifically, we propose a variant of the classical \emph{asynchronous} consensus protocol and show that we can accelerate the convergence
considerably by a factor depending on the underlying network. The main idea is to use the distributed algorithm we developed for effective resistances to design a {weight matrix} which can help pass the information among neighbors more {\it effectively} -- an alternative approach in~\cite{olshevsky2014linear} also builds on modifying the weights depending on the degree of the neighbors. 
Since the consensus iterations are the building block of many existing core distributed optimization algorithms such as the distributed subgradient, distributed proximal gradient and ADMM methods; we believe that our method and framework have far-reaching potential for accelerating many other distributed algorithms in addition to consensus algorithms, and this will be the subject of future work.}\vspace*{-0.5mm}

{\textbf{Outline.} In Section \ref{sec-method}, we introduce our algorithm for computing effective resistances. In Section \ref{sec-numerical}, we provide numerical results; 
finally, in Section \ref{sec-future}, we give a summary of our results and discuss future work.}\vspace*{-0.5mm}

{\textbf{Notation.} Let $d_i\triangleq |\cN_i|$ denote the degree of $i\in\cN$, and $m\triangleq|\mathcal{E}|$. Throughout the paper, $\cL\in\reals^{|\cN|\times |\cN|}$ denotes the weighted Laplacian of $\cG$, i.e., $\L_{ii}=\sum_{j\in\cN_i}w_{ij}$, $\cL_{ij}=-w_{ij}$ if $j\in\cN_i$, and equal to $0$ otherwise. 
The set $\sym^n$ denotes the set of $n\times n$ real symmetric matrices. We use the notation $Z = [z_i]_{i=1}^n$ where $z_i$'s are either the columns or rows of the matrix $Z$ depending on the context. $\ones$ is the column vector with all entries equal to 1, and $\id$ is the identity matrix.}
\vspace*{-3mm}
%

\vspace{-2mm}
\section{Methodology}\label{sec-method}
\vspace*{-3mm}
Clearly, $\cL$ is symmetric and positive semidefinite; and since $\cG$ is connected, the nullspace of $\cL$ is spanned by $\ones$.
In particular, consider the eigenvalue decomposition $\cL=\sum_{i=1}^n\lambda_i u_i u_i^\top$; we have $0=\lambda_1<\lambda_2\leq \ldots\leq \lambda_n$ and $u_1=\tfrac{1}{\sqrt{n}}\ones$. Recall that we would like to compute $\cL^\dag=\sum_{i=2}^n\tfrac{1}{\lambda_i}u_i u_i^\top$ in a decentralized way. First, we are going to describe a naive way to solve this problem which would converge with a linear rate, but require storing and communicating $n\times n$ matrices among the neighboring nodes. Next, we discuss that $\cL^\dag$ can be computed in a distributed way using {the} (randomized) Kaczmarz {(RK)} method with significantly less communication burden.\vspace*{-4mm}
\subsection{A consensus-based naive method for computing $\cL^\dag$:}
\vspace*{-1mm}
Let $\theta\geq \lambda_2$ and define $\bar{\cL}\triangleq\cL+\frac{\theta}{n}\one \one^\top$, i.e., $\bar{\cL}=\theta u_1 u_1^\top+\sum_{i=2}^n\lambda_iu_iu_i^\top$; hence, $\bar{\cL}^{-1}=\cL^\dag+\frac{1}{\theta}u_1 u_1^\top$. To compute $\bar{\cL}^{-1}$, consider solving $(P):\ \min_{X\in\sym^n} f(X)\triangleq\frac{1}{2}\norm{\bar{\cL}X-\id}_F^2$. Note that $f$ is strongly convex with modulus $\lambda_2$ since $\theta\geq \lambda_2$; moreover, such $\theta$ can be chosen easily in certain cases. For instance, for unweighted $\cG$, i.e., {$w_{ij}=1$ for $(i,j)\in\cE$}, it is known that $\lambda_2\leq \min_{i\in\cN} d_i$; hence, $\theta$ could be chosen after running a min-consensus algorithm over $\cG$. To solve $(P)$ in a decentralized manner, we will exploit connectivity of $\cG$. Let $\bar{\ell}_i\in\reals^n$ be {a column vector} for $i\in\cN$ such that $\bar{\cL}={[(\bar{\ell}_i)^\top]}_{i\in\cN}$, i.e., {$({\bar{\ell}_i})^\top$} denotes the $i$-th row of $\bar{\cL}$. $(P)$ can be equivalently written as follows: \vspace*{-3mm}
{\small
\begin{align*}
(P'):~\min_{X_i\in\sym^n,~i\in\cN}\left\{\sum_{i\in\cN}\norm{X_i\bar{\ell}_i-e_i}_2^2:\ X_i=X_j\ \forall~(i,j)\in\cE\right\},\vspace*{-2mm}
\end{align*}}%
where $e_i$ denotes the $i$-th standard basis vector of $\reals^n$. Although this problem is not strongly convex in $[X_i]_{i\in\cN}$, there is a way to regularize the objective $\bar{f}([X_i]_{i\in\cN})\triangleq\sum_{i\in\cN}\norm{X_i\bar{\ell}_i-e_i}_2^2$ to make it strongly convex. Indeed, it {can be shown} that for $\alpha>0$ sufficiently large, $\bar{f}_\alpha\triangleq\bar{f}+\alpha r$ is strongly convex in $[X_i]_{i\in\cN}$, where $r([X_i]_{i\in\cN})\triangleq\sum_{(i,j)\in\cE}\norm{X_i-X_j}_F^2$; and one can equivalently consider $\min\{\bar{f}_\alpha([X_i]_{i\in\cN}):\ X_i=X_j\ (i,j)\in\cE\}$. In particular, the algorithm EXTRA in~\cite{shi2015extra} exploits a similar restricted strong convexity argument and achieves a linear convergence rate for the iterate sequence. That said, the communication overhead is the main problem with this approach of solving $(P')$. In fact, at each iteration $k$, each node $i\in\cN$ communicates its local estimate $X_i^k$ to its neighbors $\cN_i$; thus, each iteration of these consensus based methods would require $\cO(2|\cE|n^2)$ real variable communications in total, e.g., EXTRA. Next, we discuss the distributed implementation of {the RK} method to compute $\cL^\dag$, which would prove itself as a more communication efficient and practical method.\vspace*{-4mm}
\subsection{Distributed Kaczmarz method for computing $\cL^\dag$:}
Consider a \emph{consistent} system $Ax=b$, where $A=[a_i^\top]_{i=1}^m\in\reals^{m\times n}$ and $b\in\reals^m$. Suppose $A$ has no rows with all zeros, and let $x^*=\argmin\{\norm{x}_2:\ Ax=b\}$. In~\cite{gower2015stochastic}, it is shown that $x^*$ can be computed using a randomized 
Kaczmarz method. 
In particular, it follows from the results in~\cite{gower2015stochastic} that starting from $x^0\in\N(A)$, the method displayed in Algorithm~\ref{alg:RK} produces $\{x^k\}_{k\geq 1}$ such that $\mathbb{E}[\norm{x^k-x^*}_2^2]\leq \rho^k\norm{x^0-x^*}_2^2$ for $k\geq 0$ with $\rho\triangleq 1-\lambda_{\min}^+(A^\top H A)$ {where $\lambda_{\min}^+(\cdot)$ denotes the smallest positive eigenvalue} and $H=\sum_{i=1}^mp_i \frac{1}{\norm{a_i}_2^2}e_ie_i^\top$; furthermore, $1-\frac{1}{\Rank(A)}\leq \rho <1$. Note that fixing $p_i=\norm{a_i}_2^2/\norm{A}_F^2$ gives us the randomized Kaczmarz in~{\cite{strohmer2009randomized,zouzias2013randomized}}.\vspace*{-3mm}
\DontPrintSemicolon
\begin{algorithm}
\small
\textbf{Initialization:} $x^{0}\in\N(A)$ \;
\For{$k\geq 0$}{
Pick $i\in\{1,\ldots,m\}$ with probability $p_i$\;
$x^{k+1}\gets x^{k}-\frac{1}{\norm{a_i}^2}(a_i^\top x^{k}-b_i)a_i$\;
}
\caption{\small RK($\{p_i\}_{i=1}^m$) -- Randomized Kaczmarz}
\label{alg:RK}
\end{algorithm}

\vspace*{-3mm} Note $\cL\cL^\dag=\sum_{i=2}^nu_i u_i^\top$ and $\id=\sum_{i=1}^nu_i u_i^\top$; hence, $\cL\cL^\dag=\id-u_1u_1^\top=\id-\frac{1}{n}\one\one^\top$. Although {the solution set} $\{X\in\sym^n: \cL X=\id-\frac{1}{n}\one\one^\top\}$ has infinitely many elements, it is {well-known} that $\cL^\dag$ is the unique solution to
\vspace{-0.2pc}
\begin{align}
\label{eq:least-norm-system}
\cL^\dag=\argmin_{X\in\sym^n}\{\norm{X}_F:\ \cL X=B\},
\vspace*{-5mm}
\end{align}
where $B\triangleq\id-\frac{1}{n}\one\one^\top$. Let $x^l, b^l\in\reals^n$ for $l\in\cN$ {be column vectors} such that $X=[x^l]_{l\in\cN}$ and $B=[b^l]_{l\in\cN}$, i.e., $b^l=e_l-\frac{1}{n}\one$. {Note $n$ columns of $\cL^\dag$ can be computed in parallel:
\vspace{-0.3pc}
{\begin{align}
\label{eq:projection}
x^l_*\triangleq\argmin_{x\in\reals^n}\{\norm{x}_2:\ \cL x= b^l\},\quad l\in\cN, \vspace*{-8mm}
\end{align}}
i.e.}, $\cL^\dag=[x^l_*]_{l\in\cN}$. {
Since $\cL^\dag\one=\zero$, 
$x^n_*=-\sum_{l=1}^{n-1}x^l_*$. Thus, one does not need to solve for all $l\in\cN$; it suffices to compute $\{x^l_*\}_{l\in\cN\setminus\{n\}}$ and calculate $x^n_*$ from these.}

Let $\{x^{l,k}\}_{k\geq 1}$ be the sequence generated when RK implemented on \eqref{eq:projection} for $l\in\cN\setminus\{n\}$. In Algorithm~\ref{alg:step}, 
we summarized the distributed nature of RK steps assuming that each $i\in\cN$ has an exponential clock with rate $r_i>0$, and when its clock ticks, the node $i$ wakes up {and} communicates with its neighbors $j\in\cN_i$ on $\cG$. More precisely, consider the resulting superposition of these point processes, and let $\{t_k\}_{k\in\integers_+}$ be the times such that one of the clocks ticks; hence, for all $k\geq 0$, the node that wakes up at time $t_k$ is node $i$ with probability $p_i=r_i/\sum_{j\in\cN}r_i$, i.e., $\{t_k\}_{k\geq 0}$ denotes the arrival times of a Poisson process with rate $\sum_{j\in\cN}r_i$. \vspace*{-2mm}
\DontPrintSemicolon
\vspace{-2mm}
\begin{algorithm}
{\small
\textbf{Initialization:} $x_i^{l,0}\gets 0$ for $l\in\cN\setminus\{n\}$ and $i\in\cN$\;
\For{$k\geq 0$}{
At time $t_k$, $i\in\cN$ wakes up w.p. $p_i=\frac{r_i}{\sum_{j\in\cN}r_i}$\;
\For{$l\in\cN\setminus\{n\}$}{
Node $i$ requests and receives $x_j^{l,k}$ from $j\in\cN_i$\;
Node $i$ computes and sends $q_i^{l,k}$ to all $j\in\cN_i$
$q_i^{l,k}=\frac{1}{\sum_{j\in\cN_i\cup\{i\}}\cL_{ij}^2}(\sum_{j\in\cN_i\cup\{i\}}\cL_{ij}x^{l,k}_j-b^l_i)$\;
Each $j\in\cN_i\cup\{i\}$ updates $x_j^{l,k+1}\gets x_j^{l,k}-\cL_{ij}q_i^{l,k}$
}}}
\caption{\small D-RK$(\{r_i\}_{i\in\cN})$ -- Decentralized RK}
\label{alg:step}
\end{algorithm}

\vspace*{-2mm} For $k\geq 0$, let $X^k\triangleq [x^{l,k}]_{l\in\cN}$ be the concatenation of D-RK sequence, where $x^{n,k}\triangleq -\sum_{l=1}^{n-1}x^{l,k}$, and define $S=\diag(s)$ such that $s_i\triangleq\sum_{j\in\cN_i\cup\{i\}}\cL_{ij}^2$ for $i\in\cN$. According to \cite{gower2015stochastic,strohmer2009randomized}, for $r_i=s_i$, we get $H=\frac{1}{\norm{\cL}_F^2}\id$, and this implies linear convergence of $\{X^k\}_{k\geq 0}$ to $\cL^\dag$ with rate $\rho=1-\left(\frac{\lambda_{\min}^+(\cL)}{\norm{\cL}_F}\right)^2$, i.e., $\mathbb{E}[\norm{X^k-\cL^\dag}_F^2]\leq \rho^k\norm{\cL^\dag}_F^2$ for $k\geq 0$. Moreover, for each $i\in\cN$, when node $i$ wakes up, D-RK requires $2d_i (n-1)$ communications -- each communication $i$ sends/receives a real variable to/from a neighboring node in $\cN_i$; hence, at each iteration, i.e., at each time a node wakes up, the expected number of communication per iteration is $N=\sum_{i\in\cN}2 p_i d_i (n-1)\leq 2d_{\max}(n-1)$. In particular, for unweighted graphs, i.e., $w_{ij}=1$ for $(i,j)\in\cE$, we have $p_i=\frac{d_i(d_i+1)}{2m+\sum_{j\in\cN}d_j^2}$ for $i\in\cN$.

Next, instead of \eqref{eq:least-norm-system}, consider implementing D-RK on a normalized system $S^{-\sfrac{1}{2}}\cL X = S^{-\sfrac{1}{2}} B$ to obtain better convergence rate in practice -- $i$-th equation in this normalized system can be computed locally at $i\in\cN$. For this system, where all the rows have unit norm, one can set $r_i=r$ for some $r>0$ for all $i\in\cN$ -- hence, nodes wake up with uniform probability, i.e., $p_i=\frac{1}{n}$ for $i\in\cN$; for this choice of equal clock rates, $H=\frac{1}{n}\id$ and $\{X^k\}_k$ converges linearly to $\cL^\dag$ with rate ${\rho_S}\triangleq 1-\frac{1}{n}\lambda_{\min}^{+}(\cL S^{-1}\cL)$. Moreover, the expected number of communication per iteration is $N=4m\frac{n-1}{n}\leq 4m$. In all 
experiments on small world random networks -- see the definition in the numerical section, D-RK implemented on the normalized system worked much better than directly implementing it on \eqref{eq:least-norm-system} {(see Fig.~\ref{fig:pseudo-Laplacian})}. {We conjecture that for certain family of random graphs,\vspace*{-2mm}
{\small
\begin{equation}\frac{1}{n}\lambda_{\min}^+ (\cL S^{-1} \cL) \geq \left(\frac{\lambda_{\min}^+ (\cL)}{\| \cL\|_F}\right)^2
\label{ineq-lambda-min-plus}
\vspace*{-2mm}
\end{equation}}%
holds with high probability {which would directly imply that $\rho_S \leq \rho$, i.e., D-RK on the normalized system would be faster.\vspace*{-3mm}	
}}
\vspace*{-2mm} 
\vspace*{-2mm}
\section{Numerical Experiments}\label{sec-numerical}
\vspace*{-3mm}
In this chapter, first we provide numerical experiments to show that $\{R_{ij}\}_{(i,j)\in\cE}$ can be computed very efficiently in a decentralized fashion, and second, we demonstrate the benefits of using effective resistances in consensus algorithms.
\vspace*{-9mm}
\subsection{Decentralized computation of $\cL^\dag$} 
\begin{figure}[h!]
\begin{center}
	\includegraphics[width=0.45 \linewidth]{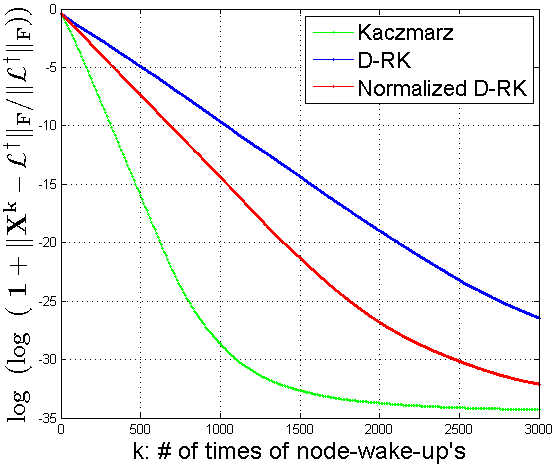}
    \includegraphics[width=0.45 \linewidth]{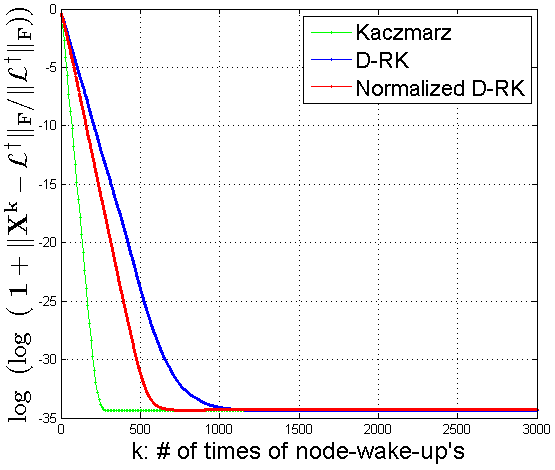}\\
    \includegraphics[width=0.45 \linewidth]{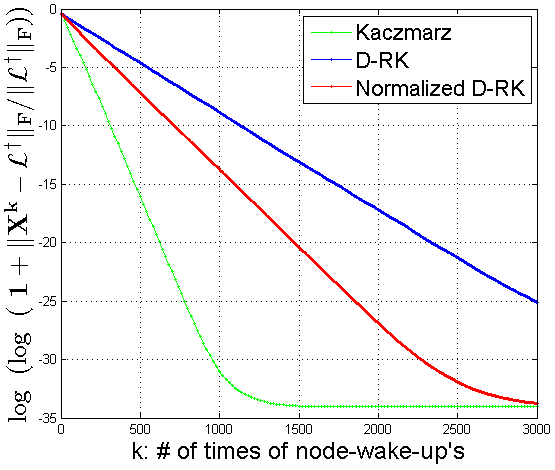}
    \includegraphics[width=0.45 \linewidth]{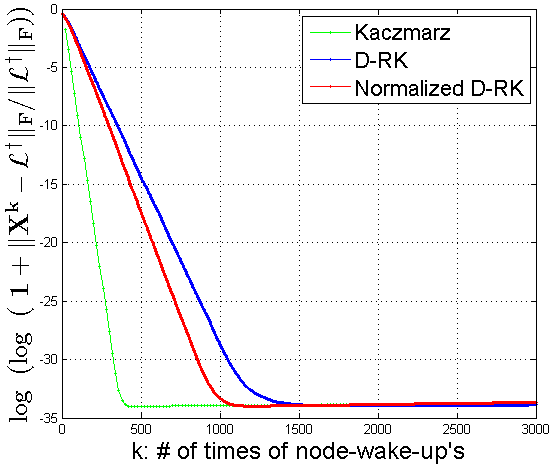}
\caption{\label{fig:pseudo-Laplacian}{\small Performance of D-RK and normalized D-RK on small-world $\cG$: \textbf{top, left:} $(n,m)=(10,18)$, \textbf{top, right:} $(n,m)=(10,36)$, \textbf{bottom, left:} $(n,m)=(20,76)$, \textbf{top, right:} $(n,m)=(20,152)$.}}
\end{center}
\vspace*{-8mm}
\end{figure}

\vspace*{-1mm}
We tested D-RK and its normalized version on unweighted small-world type communication networks, and we compared these randomized methods with deterministic (cyclic) Kaczmarz method. Given positive integers $n,m$ such that $m\geq n$, let $E\in\sym^n$ denote the adjacency matrix of the small-wold network parameterized by $(n,m)$ such that $E_{i,i+1}=1$ for $i=1,\ldots,n-1$ and $E_{1,n}=1$, and the other $m-n$ entries are chosen uniformly at random among the remaining upper diagonal elements of $E$ and set to $1$. 
We considered $n\in\{10,20\}$ and for each $n$, we chose $m$ such that the edge density, $2m/(n^2-n)$, is $0.4$ or $0.8$. For each scenario, we plot the average of $\log\log(1+\norm{X^k-\cL^\dag}_F/\norm{\cL^\dag}_F)$ {over 100 sample paths versus $k$. The results show that the randomized algorithms are slower than their deterministic counterpart; this is the price to pay for asynchronous computations. D-RK applied to the normalized system was also faster than the standard D-RK, i.e., numerically we see $\rho_S < \rho$ as suggested by the inequality \eqref{ineq-lambda-min-plus}. We also observed finite convergence on every sample path numerically -- the finite number of iterations required for convergence depended on the sample path chosen; hence, averaging iterates over sample paths led to the smooth curves reported in Fig.~\ref{fig:pseudo-Laplacian}.} \vspace*{-3mm}
\subsection{Consensus exploiting effective resistances}
Let $y^0 \in \mathbb{R}^n$ be a vector such that the $i$-th component 
represents the initial value at node $i$, and let $\bar{y} \triangleq \sum_{i=1}^n y_i^0/n$ be the average. 
In consensus algorithms, the aim is to compute $\bar{y}$ at each node in a distributed manner. As in Section~\ref{sec-method}, we assume that each $i\in\cN$ has an exponential clock with rate $r_i>0$; however, now, we assume that when its clock ticks at time $t_k$, the node $i$ wakes up and picks \emph{one} of its neighbors $j\in\cN_i$ with probability $p_{ij}\in(0,1)$, i.e., $\sum_{j\in\cN_i}p_{ij}=1$. Next, nodes $i$ and $j$ exchange their local variables $y_i^k$ and $y_j^k$. 
We assume that each node $i\in\cN$ knows $\{R_{ij}\}_{j\in\cN}$. We will be comparing two different consensus protocols, where in both protocols nodes operate as in Algorithm~\ref{alg:gossip} but with different $\{p_i\}_{i\in\cN}$ and $\{p_{ij}\}_{j\in\cN_i}$ for $i\in\cN$. \vspace*{-2mm}
\DontPrintSemicolon
\begin{algorithm}
\small
\textbf{Initialization:} $y^{0}=[y_1^{0}, y_2^{0}, \dots, y_n^{0}]^\top \in\mathbb{R}^n$ \;
\For{$k\geq 0$}{
At time $t_k$, $i\in\cN$ wakes up w.p. $p_i$\;
Picks $j\in\cN_i$ randomly w.p. $p_{ij}$\;
$y_i^{k+1}\gets \frac{y_i^{k}+y_j^{k}}{2}$,\quad $y_j^{k+1}\gets \frac{y_i^{k}+y_j^{k}}{2}$.\;
}
\vspace*{-2mm}
\caption{\small Randomized Gossiping}
\label{alg:gossip}
\end{algorithm}

\vspace*{-4mm}
\emph{Classic Randomized Gossiping:} At each iteration $k$, each edge $(i,j)\in \cE$ has equal probability of being activated. If an edge $(i,j)$ is activated at iteration $k$ the nodes take average of their decision variables $y_i^k$ and $y_j^k$. This algorithm admits an asynchronous implementation -- see, e.g., \cite{boyd2006randomized}. 
    In our node-wake-up based asynchronous setting, the same behavior can be achieved if each node $i$ wakes up with equal probability $p_i=\frac{1}{n}$, i.e., using uniform clock rates $r_i=r>0$ for $i\in\cN$, and node $i$ picks $(i,j)$ w.p. $p_{ij}=\frac{1}{d_i}$ for all $j\in\cN_i$.

\emph{Randomized Gossiping with Effective Resistances:} This algorithm is similar to classical randomized gossiping, with the only difference that edges are sampled with non-uniform probabilities
proportional to effective resistances $\{R_{ij}\}_{(i,j)\in\cN}$. In our node-wake-up based asynchronous setting, the same behavior can be achieved if each node $i$ wakes up with probability $p_i=\frac{\sum_{j\in\cN_i}R_{ij}}{2\sum_{(i,j)\in\cE}R_{ij}}$, i.e., 
setting clock rate $r_i=\sum_{j\in\cN_i}R_{ij}$ for $i\in\cN$, and node $i$ picks $(i,j)$ w.p. $p_{ij}=\frac{R_{ij}}{\sum_{j\in\cN_i}R_{ij}}$ for all $j\in\cN_i$.
\begin{figure}[h!]
\begin{center}
	\includegraphics[width=0.47 \linewidth]{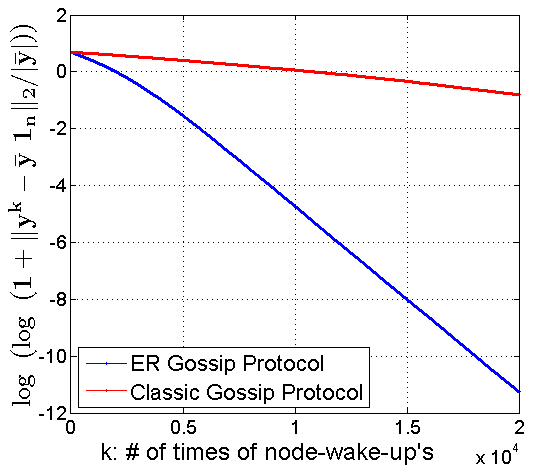}
    \includegraphics[width=0.47 \linewidth]{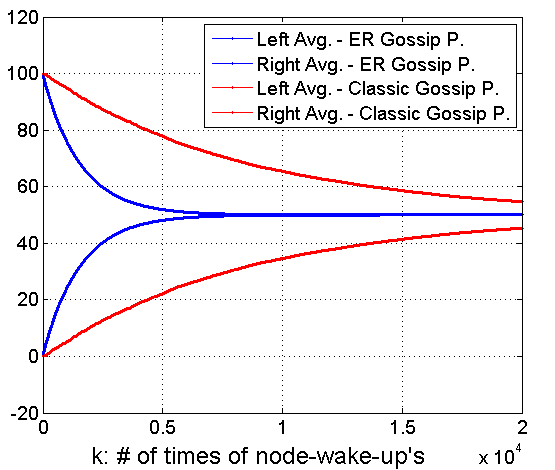}
\caption{\label{fig:gossip}{\small Performance of classic vs effective resistance based gossiping on barbell $K_{20}-K_{20}$: \textbf{left:} Relative error vs $k$, \textbf{right:} Average of left and right lobes vs $k$ for both protocols.}}
\end{center}
\vspace*{-7mm}
\end{figure}

\begin{figure}[h!]
\begin{center}
	\includegraphics[width=0.45 \linewidth]{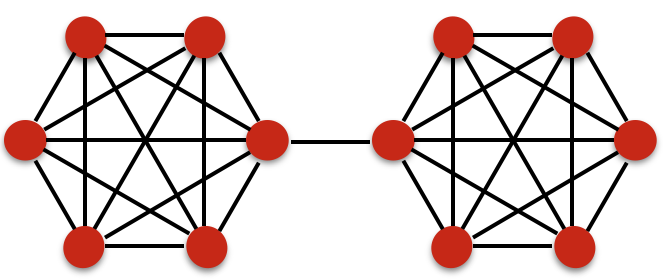}
    \caption{\label{fig-barbell}\small Barbell graph $K_n - K_n$ with $12$ nodes}
\end{center}
\vspace*{-10mm}
\end{figure}

We compare the performance of both protocols over an unweighted barbell graph $K_n - K_n$ with $2n$ nodes. Such a graph is illustrated in Fig.~\ref{fig-barbell}. In our experiment, we set $n=20$. Let $\cN_R=\{1,\ldots,20\}$ and $\cN_L=\{21,\ldots,40\}$ represent the node sets in right and left lobes {(the subgraph of $K_n$ on the right and left)} of the barbell graph. To initialize $y^0$, we sample $y_i^0$ from $\cN(100,1)$ for $i\in\cN_L$ and $y_i^0$ from $\cN(0,1)$ for $i\in\cN_R$ -- this way both lobes have significantly different local means. On the left of Fig.~\ref{fig-barbell}, we plot $\log\log(1+\norm{y^k-\bar{y}\one}_2/|\bar{y}|)$; and on the right, we plot $\tfrac{1}{20}\sum_{i\in\cN_L}y_i^k$ and $\tfrac{1}{20}\sum_{i\in\cN_R}y_i^k$ vs $k$ for both protocols. {The results show that randomized gossiping with effective resistances is much faster.}\vspace*{-3mm}
\section{Conclusions and Future Work}\label{sec-future}
\vspace{-2mm}
{In this work, we developed a distributed algorithm for computing effective resistances 
over an undirected graph $\cG$. Our method builds on an efficient, distributed and asynchronous implementation of the Kaczmarz method for solving linear Laplacian systems 
$\cL x=b$. 
We also presented an application of our algorithm to the consensus problem. }

{As part of our future work, we will investigate the finite convergence properties of 
this, suggested by the experiments. We will also study the inequality \eqref{ineq-lambda-min-plus} further
which was satisfied for a wide class of random graph models in our tests. Finally, we will investigate the applications of effective resistances to a wide class of distributed optimization algorithms which contain consensus-like iterations including distributed proximal-gradient algorithm~(DPGA) and ADMM. In particular, one could design the communication matrix $W$ for the {DPGA-W} method in~\cite{AybatMa15} using effective resistances by setting $W_{ij}=-R_{ij}$ for $(i,j) \in \cE$ and $W_{ii}=-\sum_{j\in \cN_i} R_{ij}$.
Similarly, it would be interesting to design the communication matrix in ADMM~\cite{BoydADMMBook} using effective resistances for improving its performance over for optimization problems defined over ill-conditioned graphs.}

\newpage
\bibliographystyle{IEEEbib}
\bibliography{refs_Mert,refs_Aybat}

\end{document}